\documentclass[10pt]{amsart}

\usepackage{amsmath,amsfonts, latexsym,graphicx, amssymb, amscd}

\newcommand{\U}{{\mathcal U}}
\newcommand{\0}{{\mathbf 0}}
\newcommand{\C}{{\mathbb C}}
\newcommand{\Z}{{\mathbb Z}}

\newcommand{\Pdot}{\mathbf P^\bullet}

\newtheorem{theorem}{Theorem}[section]

\newtheorem{lemma-definition}[theorem]{Lemma and Definition}
\theoremstyle{definition}

\newtheorem{corollary}[theorem]{Corollary}

\theoremstyle{remark}

\numberwithin{equation}{section}

\def\H{\mathbb H}
\def\C{\mathbb C}

\def\Z{\mathbb Z}

\begin{document}
\title {A remark on vanishing cycles with two strata}
\thanks{These notes were written with the help of the Research fund of Northeastern University}
\author{L\^e D\~ung Tr\'ang}
\address{CMI, Universit\'e de Provence\\
F-13453 Marseille cedex 13, France}
 \email{ledt@ictp.it}
\author{David B. Massey} 
\address{Department of Mathematics\\
Northeastern University, Boston, MA 02115, USA}
\email{d.massey@neu,edu}

\begin{abstract}
Suppose that the critical locus $\Sigma$ of a complex analytic function $f$ on affine space is, itself, a space with an isolated singular point at the origin $\0$, and that the Milnor number of $f$ restricted to normal slices of $\Sigma-\{\0\}$ is constant. Then, the general theory of perverse sheaves puts severe restrictions on the cohomology of the Milnor fiber of $f$ at $\0$, and even more surprising restrictions on the cohomology of the Milnor fiber of generic hyperplane slices. 
\end{abstract}

\maketitle

\section{Settings}

Let $\U$ be an open neighborhood of the origin in $\C^{n+1}$, and $f:(\U,\0)\rightarrow(\C,0)$ be a complex analytic function.
Let us call $(X,0)$ the germ of the complex analytic hypersurface defined by this function.
Suppose that the critical locus $\Sigma$ of $f$ is, itself, smooth outside of $\0$, and let $s:=\dim \Sigma$.

Throughout this paper, we make the following assumptions:

\begin{enumerate}

\item $s\geq 3$.
\item The Milnor number of a transverse slice of codimension
$s$ of the hypersurface $f^{-1}(0)$
is constant along the singular set $\Sigma$ of $X$ outside of $\0$, and equal to $\mu$.

\item The intersection of
$\Sigma$ with a sufficiently small sphere $S_\varepsilon$ centered
at $\0$ is $(s-2)$-connected.

\end{enumerate}

\medskip

\noindent Under these hypotheses, we have:

\begin{theorem}
The Milnor fiber $F_\0$ of $f$ at $\0$ can have non-zero  cohomology only in degrees $0$, $n-s$, $n-1$ and $n$.
\end{theorem}

\begin{corollary}
Suppose that $s\geq 4$ and, for a generic hyperplane $H$, the real link 
$S_\varepsilon\cap\Sigma\cap H$ of $\Sigma\cap H$ at $\0$ is $(s-3)$-connected. Then, the Milnor fiber $F_{H}$ of $f_{|_H}$ at $\0$ can have non-zero cohomology only in degrees $0$, $n-s$ and $n-1$.
\end{corollary}

\section{An exact sequence}

Let ${\Z^\bullet_{\U}}$ be the constant sheaf on $\U$ with stalks isomorphic to
the ring of integers $\Z$. If $\phi_f$ is the functor of vanishing cycles
of $f$, we know (see, e.g., \cite{D}, Theorem 5.2.21) that the complex $\phi_f[-1]{\Z^\bullet_{\U}}[n+1]$ 
is a perverse sheaf 
(see, e.g., \cite{BBD} p. 9) on $f^{-1}(0)$.
Let $\Pdot$ denote the restriction of this sheaf to its support $\Sigma$, which is the set of critical points of
$f$ inside $f^{-1}(0)$. 

We know that, for all $x\in\Sigma$, we have
$$\H^{-k}({\mathbb B}(x)\cap \Sigma; \Pdot) \ \cong  \ H^{-k}(\Pdot)_x \ \cong \ \widetilde H^{n-k}(F_x;\Z),
$$
where $F_x$ is the Milnor fiber of $f$ at $x$ and ${\mathbb B}(x)$
is a sufficiently small ball (open or closed, with non-zero radius) of ${\mathbb C}^{n+1}$ centered at $x$. Let ${\mathbb B}^*(x)={\mathbb B}(x)-\{x\}$.

Then, we have the exact sequence in hypercohomology:

\vskip.1in
$
\rightarrow {\mathbb H}^{-k}({\mathbb B}(x)\cap\Sigma, {\mathbb B}^*(x)\cap\Sigma;\Pdot)\rightarrow
{\mathbb H}^{-k}({\mathbb B}(x)\cap\Sigma;\Pdot)\hfill $
$\hbox{}\hfill \rightarrow 
{\mathbb H}^{-k}({\mathbb B}^*(x)\cap\Sigma;\Pdot)\rightarrow
{\mathbb H}^{-k+1}({\mathbb B}(x)\cap\Sigma,{\mathbb B}^*(x)\cap\Sigma;\Pdot)\rightarrow$

\vskip.1in
Since $\Pdot$ is perverse, using the cosupport condition (see e.g. \cite{BBD} p. 9):
$${\mathbb H}^{-k+1}({\mathbb B}(x)\cap\Sigma,{\mathbb B}^*(x)\cap\Sigma;\Pdot)=0$$
for $-k+1<0$. Therefore,
$$\widetilde H^{n-k}(F_x;\Z) \ \cong \ {\mathbb H}^{-k}({\mathbb B}(x)\cap\Sigma;\Pdot) \ \cong \ 
{\mathbb H}^{-k}({\mathbb B}^*(x)\cap\Sigma;\Pdot)$$
for $-k+1<0$.

\section{Topological hypothesis}

Now let us suppose that the real link $S_\varepsilon\cap\Sigma$ of the critical locus $\Sigma$ at $\0$ is $(s-2)$-connected; in particular, as $s\geq 3$, $S_\varepsilon\cap\Sigma$ is simply-connected.

 Our assumption on the constancy of the Milnor number of $f$, restricted to a normal slice to $\Sigma$, is equivalent to saying that $\Pdot_{|_{\Sigma-\{\0\}}}$ is locally constant, with stalk cohomology $\Z^\mu$ concentrated in degree $-s$. 
  As $\mathbb B^*(\0)\cap\Sigma$ is homotopy-equivalent to $S_\varepsilon\cap \Sigma$, which is simply-connected, it follows that $\Pdot_{|_{B^*(\0)\cap\Sigma}}$ is isomorphic to the shifted constant sheaf $\left(\Z^\mu\right)^\bullet_{B^*(\0)\cap\Sigma}[s]$.

This implies that
$${\mathbb H}^{-k}({\mathbb B}^*(\0)\cap\Sigma;\Pdot) \ \cong  \ 
{H}^{-k+s}({\mathbb B}^*(0)\cap\Sigma;{\mathbb Z}^\mu) \ \cong  \ 
{H}^{-k+s}(S_\varepsilon\cap\Sigma;{\mathbb Z}^\mu).$$

Thus, as $S_\varepsilon\cap\Sigma$ is $(s-2)$-connected, we have:
$${\mathbb H}^{-s}({\mathbb B}^*(\0)\cap\Sigma;\Pdot) \ \cong \ 
{H}^{0}(S_\varepsilon\cap\Sigma;{\mathbb Z}^\mu) \ \cong \ {\mathbb Z}^\mu,$$
and, if $2\leq k\leq s-1$:
$${\mathbb H}^{-k}({\mathbb B}^*(\0)\cap\Sigma;\Pdot) \ \cong \ 
{H}^{s-k}(S_\varepsilon\cap\Sigma;{\mathbb Z}^\mu)=0.$$

\section{Proofs}

Combining the results from the previous two sections, we find that, if the real link of the critical locus $\Sigma$ at $\0$ is $(s-2)$-connected
and  $s\geq 3$, then we have for the Milnor fiber $F$ of $f$ at $\0$:
$$\widetilde H^{n-s}(F; {\mathbb Z}) \ \cong  \ H^0(S_\varepsilon\cap\Sigma;{\mathbb Z}^\mu)
 \ \cong \  {\mathbb Z}^\mu;$$
$$\widetilde H^{k}(F;{\mathbb Z})=0 \hbox{, if $k\neq n-1,n$ }.$$
This proves the theorem.

Suppose now that, in addition to our other hypotheses, $s\geq 4$ and, for generic hyperplanes $H$, $S_\varepsilon\cap\Sigma\cap H$ is $(s-3)$-connected. Then, $f_{|_H}$ satisfies the hypotheses of the theorem, except that $n$ is replaced by $n-1$ and $s$ is replaced by $s-1$.  Thus, for the Milnor fiber $F_{H}$:

$$\widetilde H^{n-s}(F_{H}; {\mathbb Z}) \ 
 \ \cong \  {\mathbb Z}^\mu;$$
$$\widetilde H^{k}(F_{H};{\mathbb Z})=0 \hbox{, if $k\neq n-2,n-1$ }.$$

However, by the main result of \cite{L}, the Milnor fiber $F$ is obtained from the Milnor fiber $F_H$ by attaching cells in dimension $n$. Hence, $\widetilde H^{n-2}(F_{H};{\mathbb Z})\cong \widetilde H^{n-2}(F;{\mathbb Z})$, which we know is $0$. This proves the corollary.

\section{When the critical locus is an ICIS}

Assume that the critical locus $\Sigma$ of $f$ is an isolated complete intersection singularity 
 (ICIS) of dimension $s\geq 4$. 
 
 For an ICIS, the real link
$S_\varepsilon\cap\Sigma$ is $(s-2)$-connected. In addition, for a generic hyperplane $H$, the critical locus of $f_{|_H}$, which equals $\Sigma\cap H$, will also be an ICIS, but now of dimension $s-1$. Thus, $S_\varepsilon\cap \Sigma\cap H$ is $((s-1)-2)$-connected. Therefore, we are in the situation 
that we have considered above.

In his preprint \cite{S} M. Shubladze asserts that if the singular locus $\Sigma$ of $f$ is a complete
intersection with isolated singularity at $\0$ of dimension $\geq 3$ and the Milnor number
for transverse sections is $1$ along $\Sigma\setminus\{\0\}$, the Milnor number of $f$ at $0$
has cohomology possibly $\neq 0$ only in dimensions $0$, $n-s$ and $n$.

The results above show that, under the hypothesis of M. Shubladze, one obtains in a general way
that the cohomology of the Milnor fiber of $f$ at $\0$ is 
possibly $\neq 0$ in dimension $0$, $n-s$, $n-1$ and $n$, and
a similar result as the one of M. Shubladze in dimension $0$, $n-s$, $n-1$ for the cohomology 
of the Milnor fiber of $f$ restricted to a 
general hyperplane section if $\dim \Sigma\geq 4$.

\smallskip
In light of our corollary, it seems reasonable to ask: Can every function such as that studied by Shubladze can be obtained as a generic hyperplane restriction of a function satisfying the same hypotheses?

\smallskip

Shubladze's proof is via deformation and does not seem to easily answer the question above.

\section{What if $S_\varepsilon\cap\Sigma$ is a homology sphere?}

One might also wonder what happens if the real link of  $\Sigma$ is $(s-1)$-connected. This would, in fact, imply that $S_\varepsilon\cap\Sigma$ is a homology sphere. In this case, our earlier exact sequence immediately yields that $\widetilde H^{n-1}(F;\Z)=0$.

\smallskip

A special case of $S_\varepsilon\cap\Sigma$ being a homology sphere would occur if $\Sigma$ 
were smooth. However, in this case, Proposition 1.31 of \cite{M} implies that the Milnor 
number cannot change at $\0$, i.e., we have a smooth $\mu$-constant family, and so 
the non-zero cohomology of $F$ occurs only in degrees $0$ and $n-s$.


\begin{thebibliography}{A}
\bibitem{BBD} A. A. Beilinson, J. Bernstein, P. Deligne  Faisceaux pervers, AstŽrisque {bf 100}
(1982), Soci\'et\'e Math\'ematique de France, Paris. 
\bibitem{D} A. Dimca, Sheaves in Topology, Springer-Verlag, 2004.
\bibitem{L} L\^e D\~ung Tr\'ang, Calcul du nombre de cycles \'evanouissants d'une hypersurface complexe, Ann. Inst. Fourier (Grenoble) {\bf 23},  261--270, 1973.
\bibitem{M} D. Massey, L\^e Cycles and Hypersurface Singularities, Lecture Notes in Mathematics {\bf 1615}, Springer-Verlag, 1995.
\bibitem{S} M. Shubladze, Singularities with critical locus an complete intersection 
and transversal type $A_1$, Preprint 2010-16, Max Planck Institute for Mathematics preprints, http://www.mpim-bonn.mpg.de/preblob/4141
\end{thebibliography}
 \end{document}